\documentclass{article}
\usepackage{latexsym}
\usepackage{amsmath,amssymb}
\usepackage[dvips]{graphicx, color}

\begin{document}
\title{\bf\Large{{COUNTING\ \ \  BUNDLES} \qquad}}
\author{\bf{\large{Lin WENG}}}
\date{}
\maketitle

\section{Zetas for $(G,P)/\mathbb F_q(X)$}

Let $X$ be an irreducible reduced regular projective curve defined over $\mathbb F_q$ with $F$ the field of rational functions. Let $G$ be a reductive group of rank $r$ with $B$ a fixed Borel subgroup, both defined over $F$.
As usual, $\Delta$ stands for the corresponding collection of simple roots; $W$ the associatd Weyl group;  for a 
root $\alpha$, $\alpha^\vee$ the corresponding coroot; and  $\rho:=\frac{1}{2}\sum_{\alpha>0}\alpha$, the Weyl vector\,$\dots$
\vskip 0.20cm
\noindent
{\bf Definition 1} The {\it period for $G$ over $F$} is defined by, for $\mathrm{Re}\,\lambda\in\mathcal C^+$,
$$\boxed{\omega_{F}^G(\lambda):=\sum_{w\in W}\Bigg(\frac{1}{\prod_{\alpha\in\Delta}(1-q^{-\langle w\lambda-\rho,\alpha^\vee\rangle})}\cdot\prod_{\alpha>0,w\alpha<0}\frac{\widehat{\zeta}_{F}(\langle\lambda,\alpha^\vee\rangle)}
{\widehat{\zeta}_{F}(\langle\lambda,\alpha^\vee\rangle+1)}\Bigg)}$$ where $\mathcal C^+$ denotes the so-called positive chamber of $\frak a_0$, the space of 
characters associated to $(G,B)$, and 
$\widehat{\zeta}_{F}(s)$ the complete Artin zeta function of $X/{\mathbb F}_q$.

Corresponding to a fixed maximal standard parabolic subgroup $P$ is a simple root $\alpha_P\in\Delta$. 
Write $\Delta\backslash\{\alpha_P\}=:\{\beta_{1,P}, \beta_{2,P},\dots, \beta_{r-1,P}\}.$ 
\vskip 0.20cm
\noindent
{\bf Definition 2} The {\it period for $(G,P)$ over $F$} is defined by, for $\mathrm{Re}\lambda_P\gg 0$,
$$\boxed{\omega_{F}^{G/P}(\lambda_P):=
\mathrm{Res}_{\langle \lambda-\rho,\beta_{r(G)-1,P}^\vee\rangle=0}\cdots
\mathrm{Res}_{\langle \lambda-\rho,\beta_{2,P}^\vee\rangle=0}\mathrm{Res}_{\langle \lambda-\rho,\beta_{1,P}^\vee\rangle=0}\,\Big(\omega_{F}^G(\lambda)\Big)}$$
Here, starting from $r$-variable $\lambda\in\frak a_{0,\mathbb C}^*$, after taking residues along with 
 $(r-1)$ (independent) singular hyperplanes
$$\langle \lambda-\rho,\beta_{1,P}^\vee\rangle=0,\,
\langle \lambda-\rho,\beta_{2,P}^\vee\rangle=0,\,
\cdots,\,
\langle \lambda-\rho,\beta_{r(G)-1,P}^\vee\rangle=0,$$
we are left with only one variable, which we call $\lambda_P$.
\vskip 0.20cm
Clearly, there is a minimal integer $I(G/P)$ and finitely many factors (depending on the choice of $\lambda_P$), 
$$\widehat{\zeta}_{F}\Big(a_1^{G/P}\lambda_P+b_1^{G/P}\Big),\ \widehat{\zeta}_{F}\Big(a_2^{G/P}\lambda_P+b_2^{G/P}\Big),\ \cdots,\ 
\widehat{\zeta}_{F}\Big(a_{I(G/P)}^{G/P}\lambda_P+b_{I(G/P)}^{G/P}\Big),$$
such that {\it no $\widehat{\zeta}_{F}(a\lambda_P+b)$ factors appear in the denominators of (all terms of)
the product $\Big[\prod_{i=1}^{I(G/P)}\widehat{\zeta}_{F}\Big(a_i^{G/P}\lambda_P+b_i^{G/P}\Big)\Big]\cdot \omega_{F}^{G/P}(\lambda_P)$.}

\vskip 0.30cm
\noindent
{\bf Definition 3} (i) The {\it zeta function $\widehat{\zeta}_{F;o}^{G/P}$ for $(G,P)$ over $F$} is defined by
$$\boxed{\widehat{\zeta}_{F;o}^{G/P}\Big(s\Big)
:=\Bigg[\prod_{i=1}^{I(G/P)}\widehat{\zeta}_{F}\Big(a_i^{G/P}s+b_i^{G/P}\Big)\Bigg]\cdot \omega_{F}^{G/P}\Big(s\Big),\ \ \mathrm{Re}\,s\gg 0}
$$

\noindent
{\bf Functional Equation} {\it There exists a constant $c_{G/P}\in \mathbb Q$ such that} $$\boxed{\widehat{\zeta}_{F;o}^{G/P}\Big(-s+c_{G/P}\Big)=\widehat{\zeta}_{F;o}^{G/P}\Big(s\Big).}$$

\vskip 0.30cm
\noindent
{\bf Definition 3} (ii) The {\it zeta function $\widehat{\zeta}_{F}^{G/P}\Big(s\Big)$ for $(G,P)$ over $F$} is defned by
$$\boxed{\widehat{\zeta}_{F}^{G/P}\Big(s\Big):=\widehat{\zeta}_{F;o}^{G/P}\Big(s+\frac{c_{G/P}-1}{2}\Big)}$$
\section{Non-Abelian Zeta Functions for $\mathbb F_q(X)$}
Motivated by [W1], and for the RH, introduce a new genuine {\it pure} non-abelian  zetas for $X$ by
$$\boxed{\zeta_{X,r}(t):=\sum_{m=0}^\infty
\sum_{V\in [V]\in\mathcal M_{X,r}(d), d=rm}\frac{q^{h^0(X,V)}-1}{\#\mathrm{Aut}(V)}\cdot (q^{-s})^{d(V)},\qquad\mathrm{Re}(s)>1}$$ Here, as usual, $\mathcal M_{X,r}(d)$ denotes the moduli space of semi-stable $\mathbb F_q$-rational vector bundles of rank $r$, $[\ ]$ the Seshedri class defined using Jordan-H\"older graded bundles, and $\mathrm{Aut}(V)$ denotes the automorphism group of $V$.
Introduce also the complete zeta function for $X$ by
$$\boxed{\widehat {Z}_{X,r}(t):=\sum_{m=0}^\infty
\sum_{V\in [V]\in\mathcal M_{X,r}(d), d=rm}\frac{q^{h^0(X,V)}-1}{\#\mathrm{Aut}(V)}\cdot (t^s)^{\chi(X,V)}}$$ 
\vskip 0.30cm
\noindent
{\bf Zeta Facts} (i) ({\bf Rationality})
$$\boxed{\begin{aligned}Z_{X,r}(t)
=&\sum_{m=0}^{(g-1)-1}\alpha_{X,r}(mr)
\cdot \Big(t^{rm}+ q^{r[(g-1)-m]}\cdot t^{r[2(g-1)-m]}\Big)\\
&+\alpha_{X,r}(r(g-1))\cdot t^{r(g-1)}+\beta_{X,r}(0)\cdot \frac{(q^r-1)t^{rg}}{(1-q^rt^r)(1-t^r)}
\end{aligned}}$$
{\it with}
$$\boxed{\beta_{X,r}(0):=\sum_{V\in [V]\in\mathcal M_{X,r}(0)}\frac{1}{\#\mathrm{Aut}(V)},\ \ \ \alpha_{X,r}(d):=\sum_{V\in [V]\in\mathcal M_{X,r}(d)}\frac{q^{h^0(X,V)}-1}{\#\mathrm{Aut}(V)}}$$

\noindent
(ii) ({\bf Functional Equation}) $$\boxed{\widehat Z_{X,r}({1}/{qt})=\widehat Z_{X,r}(t)}$$
We  expect to have the following
\vskip 0.10cm
\noindent
{\bf Riemann Hypothesis}
$$\boxed{All\ zeros\ of\ the\ zeta\ function\ \widehat{\zeta}_{X,r}(s)\ lie\ on\ the\ central\ line\
\mathrm{Re}\,s=\displaystyle{\frac{1}{2}}}$$

\section{Counting Bundles}

{\color[cmyk]{0,1,1,0}Semi-stable bundles are naturally counted within the stratifications of the refined Brill-Noether loci defined using
$h^0$ and $\mathrm{Aut}$. So the Riemann Hypothesis offers us intrinsic  quantitative controls uniformly
through $\alpha$'s and $\beta$.} More generally, 
write {\it ss} for the part corresponding to {\it semi-stable principal bundles}.
\vskip 0.30cm
\noindent
{\bf Counting Conjectures} (i) ({\bf Parabolic Reduction, Stability \& the Mass})
$$\boxed{\mathrm{Vol}\Big(K_{\mathcal O}Z_{G(\mathbb A)}\backslash G^1(\mathbb A)_{\mathrm{ss}}/G(F)\Big)=\mathrm{Res}_{\lambda=\rho}\omega^G_F(\lambda)=
\mathrm{Res}_{s=1}\widehat{\zeta}_{F}^{G/P}\Big(s\Big)}$$
(ii) ({\bf Uniformity}) {\it There exist rational functions $R_{r,q}(t)$ depending on $q$ and $t$ and rational numbers $a_r, b_r$ depending only on $r$, but
independent of the curve $X$, such that} $$\boxed{\widehat{\zeta}_{F,r}(s)=R_{r,q}(q^{-a_rs-b_r})\cdot \widehat{\zeta}_{F}^{SL_r/P_{r-1,1}}\Big(a_rs+b_r\Big)}$$
\ \  Parallel structures for number fields are exposed in [W2, 3],  [Ko], [KKS], and [W5] which contains an adelic approach to Atiyah-Bott \& Witten ([AB], [Wi]) and to Kontsevich ([K]). For function fields, (i) is related to Harder-Narasimhan ([HN],  [Z], [LR]), uniformity holds for $G=SL_2$ ([W4]), for which, the RH is established in [Y], and a proof of [Ko] style for the FE of general $\widehat{\zeta}^{G/P}_F(s)$ can be obtained.

\vskip 0.40cm
\centerline{\bf REFERENCES}
\vskip 0.10cm
\noindent
[A1,2] J. Arthur, A trace formula for reductive groups. I. 
Terms associated to classes in $G({\mathbb Q})$. Duke Math. J. {\bf 45} (1978), 
no. 4, 911--952; II. Applications of a 
truncation operator. Compositio Math. {\bf 40} (1980), no. 1, 87--121.
\vskip 0.10cm
\noindent
[A3] J. Arthur, A measure on the unipotent variety, Canad. J. Math {\bf 37},
(1985) pp. 1237--1274
\vskip 0.10cm
\noindent
[AB] M. Atiyah \& R. Bott, The Yang-Mills equations over Riemann surfaces, Philos. Trans. Roy. Soc. London 308, 523-615 (1983)
\vskip 0.10cm
\noindent
[B] K. Behrend, The Lefschetz trace formula for the moduli stack of principal
bundles. PhD thesis, UC Berkley, 1990
\vskip 0.10cm
\noindent 
[HN] G. Harder \& M.S. Narasimhan, On the cohomology groups of moduli spaces of vector bundles on curves, Math. Ann. 212, 215-248 (1975)
\vskip 0.10cm
\noindent 
[JLR] H. Jacquet, E. Lapid \& J. Rogawski,  Periods of automorphic forms. 
J. Amer. Math. Soc. 12 (1999), no. 1, 173--240
\vskip 0.10cm
\noindent 
[KKS] H. Ki, Y. Komori \& M. Suzuki, On the zeros of Weng zeta functions for Chevalley groups, arXiv:1011.4583v1
\vskip 0.10cm 
\noindent  
[KW] H.H. Kim \& L. Weng, Volume of truncated fundamental domains, Proc. AMS 135, 1681-1688 (2007)
\vskip 0.10cm
\noindent
[KS] M. Kontsevich \& Y. Soibelman, Lectures on motivic Donaldson-Thomas invariants and Wall-crossing formulas, manuscripts, Dec.  2011
\vskip 0.10cm
\noindent
[Ko] Y. Komori, Functional equations for Weng's zeta functions for $(G,P)/\mathbb Q$, Amer. J. Math., to appear
\vskip 0.10cm
\noindent
[L] L. Lafforgue, {\it Chtoucas de Drinfeld et conjecture de 
Ramanujan-Petersson}. Asterisque No. 243 (1997)
\vskip 0.10cm
\noindent
[La] R. Langlands, The volume of the fundamental domain for some arithmetical 
subgroups of Chevalley groups, in {\it Algebraic Groups and Discontinuous 
Subgroups,} Proc. Sympos. Pure Math. 9, AMS (1966) pp.143--148
\vskip 0.10cm
\noindent
[LR] G. Laumon\& M. Rapoport, The Langlands lemma and the Betti numbers of stacks of G-bundles on a curve, Intern. J. Math. 7(1996), 29-45.
\vskip 0.10cm 
\noindent  
{[We]} A. Weil,  {\it Adeles and algebraic groups},  Progress in Math., 23. 1982.
\vskip 0.10cm
\noindent
[W1] L. Weng, Non-abelian zeta function for function fields, Amer. J. Math., 127 (2005), 973-1017
\vskip 0.10cm
\noindent
[W2] L. Weng, A geometric approach to $L$-functions, in {\it Conference on L-Functions}, pp. 219-370, World Sci (2007)
\vskip 0.10cm
\noindent
[W3] L. Weng, Symmetry and the Riemann Hypothesis, in {\it Algebraic and Arithmetic Structures of Moduli Spaces}, ASPM 58, 173-223 (2010)
\vskip 0.10cm
\noindent
[W4] L. Weng, Zeta Functions for Elliptic Curves I: Counting Bundles, preprint, 2012, arXiv:1202.0870
\vskip 0.10cm
\noindent
[W5] L. Weng, Parabolic Reduction, Stability \& the Mass, in preparation
\vskip 0.10cm
\noindent
[Wi] E. Witten, On quantum gauge theories in two dimensions, Commun. Math. Phys, 141 (1991) 153-209
\vskip 0.10cm
\noindent
[Y] H. Yoshida, manuscripts, Dec., 2011
\vskip 0.10cm
\noindent
[Z] D. Zagier, Elementary aspects of the Verlinde formula and the Harder-Narasimhan-Atiyah-Bott formula, in {\it Proceedings of the Hirzebruch 65 Conference on Algebraic Geometry}, 445-462 (1996)
\vskip 0.30cm
\noindent
{\bf Lin WENG}\footnote{{\it Acknowledgement.}  We would like to thank  M. Kontsevich and H. Yoshida for sharing with us their works, and C. Deninger and H. Hida for their constant encouragements.

This work is partially supported by JSPS.
} Institute for Fundamental Research, The $L$-Academy {\it and}

\noindent
Graduate School of Mathematics,  Kyushu University, 
Fukuoka 819-0395, 
Japan

\noindent
E-Mail: weng@math.kyushu-u.ac.jp
\end{document}